\documentclass[letterpaper, 11pt]{article}

\usepackage{graphics,graphicx,xcolor,float,enumerate}
\usepackage{multicol} 
\usepackage{parskip}
\usepackage{amsmath}
\usepackage{multirow}
\usepackage[utf8]{inputenc}
\usepackage{fancyhdr}
\usepackage[title]{appendix}
\usepackage{wasysym}
\usepackage{url}
\usepackage{mathtools}
\usepackage{booktabs}

\usepackage[font=footnotesize,labelfont=small]{caption}
\RequirePackage{geometry}
\newcommand{\mtwo}[4]{\left[\begin{array}{cc} #1 & #2 \\ #3 & #4 \end{array}\right]}
\definecolor{dgreen}{rgb}{0.0, 0.0, 0.0}
\definecolor{dbrown}{rgb}{0.0, 0.0, 0.0}
\definecolor{dblue}{rgb}{0.0, 0.0, 0.0}

\newcommand{\mtx}[1]{\boldsymbol{#1}}
\newcommand{\vct}[1]{\boldsymbol{#1}}

\newcommand{\ycmt}[1]{\textcolor{blue}{#1}}
\title{An alternative extended linear system for boundary value problems on 
locally perturbed geometries}
\author{Y. Zhang and A. Gillman}

\date{}

\begin{document}
\maketitle
\begin{abstract}
   This manuscript presents a new extended linear system for integral equation
   based techniques for solving 
   boundary value problems on locally perturbed geometries. 
   \textcolor{dbrown}{The new extended linear system is similar to 
   a previously presented technique for which the authors have constructed
   a fast direct solver.  The key features of the work presented in 
   this paper are that the fast direct solver is more efficient for 
   the new extended linear system and that problems involving 
  specialized quadrature for weakly singular kernels can be easily handled. } 
   Numerical results illustrate the improved performance of the fast direct 
   solver for the new extended system when compared to the fast direct solver
   for the original extended system.
\end{abstract}

\section{Introduction}
This manuscript presents an integral equation based solution technique for elliptic boundary value problems
on locally-perturbed geometries. Such problems arise in applications such as optimal shape design.
\textcolor{dbrown}{In each iteration or optimization cycle, the changes to the object shape often stay local to certain parts of the object.}
The proposed approach formulates an extended linear system that allows for the boundary value problem on the new geometry 
to be expressed in terms of a linear system on the original geometry plus a correction to account for the local perturbation.
This idea was first proposed in \cite{2009_martinsson_ACTA} and a fast direct solver was constructed for 
the resulting formulation in \cite{ZHANG2018_jcp}.
Unfortunately, the fast direct solver for the original extended system required inverting a matrix the size of 
the number of points removed from the original geometry which is expensive if the removed portion is 
large.  
Another difficulty of the
original extended system is that care is required when the technique is applied to systems discretized
using quadrature for weakly singular kernels.  
The extended linear system proposed in this manuscript overcomes these two difficulties.  Additionally, a fast direct solver for the new extended system can 
be constructed from the tools presented in \cite{ZHANG2018_jcp} but is more efficient than the original fast 
direct solver. 

This manuscript briefly reviews a boundary integral formulation for a Laplace problem with Dirichlet boundary data 
and the linear system that results from the discretization in Section \ref{sec:BIE}.  Next, the original extended
system and the new extended system are presented in Section \ref{sec:extended}. Finally numerical results illustrate
the efficiency of the fast direct solver for the new extended system in Section \ref{sec:numerics}.

\section{Boundary integral formulation}
\label{sec:BIE}

\begin{figure}[ht]
\begin{center}
\begin{tabular}{c}
 \setlength{\unitlength}{1mm}
\begin{picture}(60,30)(0,0)
\put(15,0){\includegraphics[height=30mm]{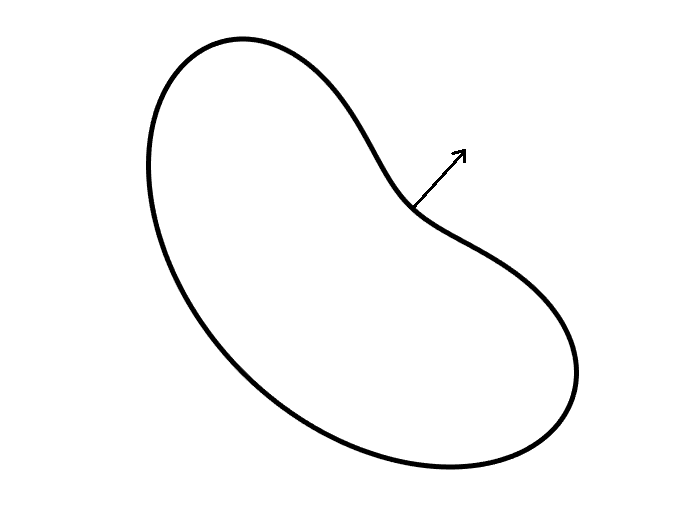}}
 \put(30,12){$\Omega$}
\put(28,29){$\Gamma$}
\put(35,16){$\vct{x}$}
\put(42,22){$n({\vct{x}})$}
\end{picture}\\
(a)
\end{tabular}
\begin{tabular}{c}
 \setlength{\unitlength}{1mm}
\begin{picture}(60,30)(0,0)
\put(15,0){\includegraphics[height=30mm]{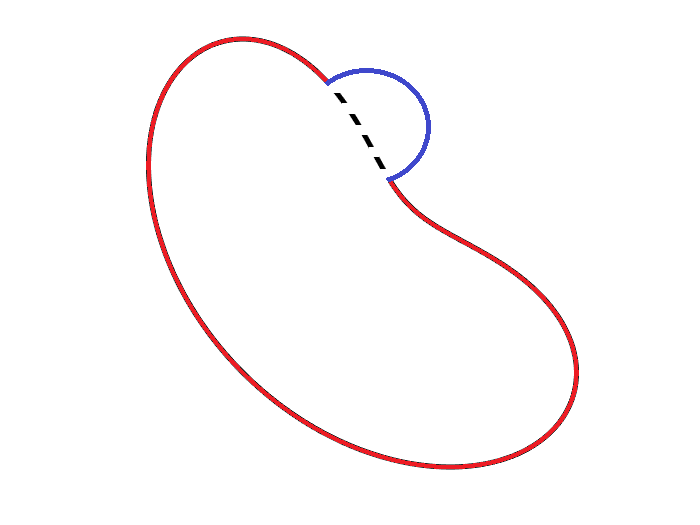}}
 \put(17,20){$\Gamma_{k}$}
\put(35,28){$\Gamma_{p}$}
\put(30,20){$\Gamma_{c}$}
\end{picture}\\
(b)
\end{tabular}
\end{center}
\caption{ (a) A sample geometry $\Omega$ with boundary $\Gamma$ and 
outward facing normal vector $n(\vct{x})$ at the point $\vct{x}\in \Gamma$.  (b) A sample locally perturbed geometry where 
the original boundary is $\Gamma_o = \Gamma_k\cup \Gamma_c$, the portion of the boundary
being removed is $\Gamma_c$, the portion of the original boundary remaining is $\Gamma_k$ 
and the newly added boundary is $\Gamma_p$.}
\label{fig:laplaceBvp}
\end{figure}
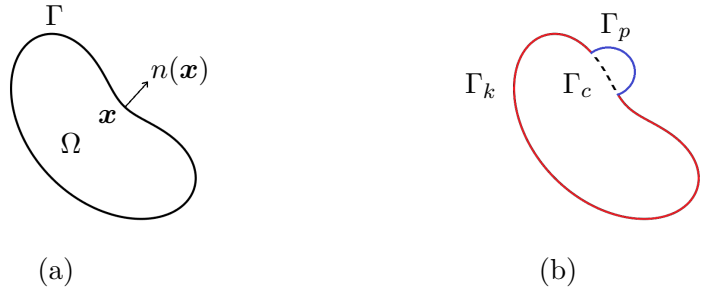

Consider the interior Laplace problem with Dirichlet boundary condition
\begin{equation}\label{eq:laplaceBvp}
\begin{array}{rlll}
 -\Delta u(\vct{x}) &=&0 &\mbox{ for }\vct{x} \in \Omega,\\
 u(\vct{x})&=&g(\vct{x})&\mbox{ for }\vct{x} \in \Gamma\\
\end{array}
\end{equation}
where $\Omega$ denotes the interior of the geometry, and $\Gamma$ denotes the boundary of $\Omega$,
as illustrated in Figure \ref{fig:laplaceBvp}(a).
Let $G(\vct{x},\vct{y}) = -\frac{1}{2\pi}\log|\vct{x}-\vct{y}|$ denote the Green's function for the Laplace operator and
$D(\vct{x},\vct{y}) = \partial_{n(\vct{y})}G(\vct{x},\vct{y})$ denote the double layer kernel where 
$n(\vct{x})$ denotes the outward facing normal vector at the point $\vct{x}\in\Gamma$.  
The solution to (\ref{eq:laplaceBvp}) can be expressed as 
\begin{equation}
 \label{equ:laplaceSoln}
u(\vct{x}) = \int_{\Gamma}  D(\vct{x},\vct{y})\sigma(\vct{y}) ds(\vct{y}) \qquad\mbox{for } \vct{x}\in \Omega,
\end{equation}
where $\sigma(\vct{x})$ is some unknown density defined only on the boundary $\Gamma$. 
Enforcing that $u(\vct{x})$ satisfies the boundary condition results in the following integral equation 
for $\sigma(\vct{x})$; 
\begin{equation}
-\frac{1}{2}\sigma(\vct{x}) +\int_{\Gamma} D(\vct{x},\vct{y})\sigma(\vct{y}) ds(\vct{y}) = g(\vct{x}).
\label{equ:laplaceBie} 
\end{equation}

Upon discretization via a Nystr\"om or boundary element method, one is left with solving a 
dense linear system 
\begin{equation}
\mtx{A}\vct{\sigma} = \vct{g},
\label{equ:discretizedBie}
\end{equation}
where $\mtx{A}$  is the discretized boundary integral operator and 
$\vct{\sigma}$ is the vector approximating $\sigma$ at the discretization points.

\section{Extended linear systems}
\label{sec:extended}
Consider a boundary value problem on a geometry with a local perturbation as illustrated in Figure \ref{fig:laplaceBvp}(b).
Let $\Gamma_o$ denote the boundary of the original geometry, $\Gamma_k$ denote the portion of the boundary that is not 
changing and $\Gamma_c$ denote the portion that is cut or removed. So $\Gamma_o = \Gamma_c\cup \Gamma_k$.  Let $\Gamma_p$ denote the new portion of the boundary.  Then the new geometry has a boundary $\Gamma_n$ defined by $\Gamma_n = \Gamma_k\cup\Gamma_p$.  

The discretized linear systems can be partitioned according to this notation.  In other words, 
the 
original system  can be expressed as
\begin{equation}\label{equ:originalLinearSystem}
\mtx{A}_{oo} \vct{\tau}_o =
\begin{bmatrix}
\mtx{A}_{kk} & \mtx{A}_{kc}\\
\mtx{A}_{ck} & \mtx{A}_{cc}\\
\end{bmatrix}
\begin{pmatrix}
\vct{\tau}_k\\
\vct{\tau}_c
\end{pmatrix}
=
\begin{pmatrix}
\vct{g}_k\\
\vct{g}_c
\end{pmatrix}
=
\vct{g}_o,
\end{equation}
and the linear system for the perturbed geometry can be expressed as
\begin{equation}\label{equ:perturbedlLinearSystem}
\mtx{A}_{nn} \vct{\sigma}_n =
\begin{bmatrix}
\mtx{A}_{kk} & \mtx{A}_{kp}\\
\mtx{A}_{pk} & \mtx{A}_{pp}\\
\end{bmatrix}
\begin{pmatrix}
\vct{\sigma}_k\\
\vct{\sigma}_p
\end{pmatrix}
=
\begin{pmatrix}
\vct{g}_k\\
\vct{g}_p
\end{pmatrix}
=
\vct{g}_n
\end{equation}
where $\vct{\tau}_k$ and $\vct{\sigma}_k$ denote the vector whose entries are the approximate solution at the discretization points on $\Gamma_k$,
$\vct{\tau}_c$ denotes the vector whose entries are the approximate solution at the discretization points on $\Gamma_c$, etc.
Likewise $\mtx{A}_{kk}$ is the submatrix of the discretized integral equation corresponding to the interaction of
$\Gamma_k$ with itself, $\mtx{A}_{kc}$ is the submatrix of the discretized integral equation corresponding to the interaction of
$\Gamma_k$ with $\Gamma_c$, etc.

\subsection{Original extended linear system}
\label{sec:orig}

The discretized problem on $\Gamma_n$ can be expressed as an extended 
linear system \cite{2009_martinsson_ACTA} by
\begin{equation}
\label{eq:orig_extended_system}
\left(\underbrace{
\mtwo{\mtx{A}_{oo}}{\mtx{0}}{\mtx{0}}{\mtx{A}_{pp}}}_{\tilde{\mtx{A}}}
+
\underbrace{\begin{bmatrix}
0 & \begin{pmatrix}-\mtx{A}_{kc}\\-\mtx{B}_{cc}\end{pmatrix} & \mtx{A}_{op}\\
\mtx{A}_{pk} & 0 & 0
\end{bmatrix}}_{\mtx{Q}_{\rm orig}}
\right)
\underbrace{
\begin{pmatrix}
\vct{\sigma}_k\\
\vct{\sigma}_c\\
\mtx{\sigma}_p
\end{pmatrix}}_{\textcolor{black}{\boldsymbol{\sigma}_{\rm ext}}}
=
\underbrace{
\begin{pmatrix}
\vct{g}_k\\ \vct{0}\\ \vct{g}_p
\end{pmatrix}
}_{\vct{g}_{\rm ext}}
\end{equation}
where $\mtx{A}_{kc}$ denotes the submatrix of $\mtx{A}_{oo}$ corresponding to
the interaction between $\Gamma_k$ and $\Gamma_c$, $\mtx{A}_{op}$ 
denotes the discretization of the double layer integral operator
on $\Gamma_p$ evaluated on $\Gamma_o$, $\mtx{A}_{pk}$ denotes 
the discretization of the double layer integral operator
on $\Gamma_k$ evaluated on $\Gamma_p$, and 
$\mtx{B}_{cc}$ denotes the sub-matrix of $\mtx{A}_{oo}$ 
corresponding to the interaction of $\Gamma_c$ with itself
but the diagonal entries are set to zero.  
The matrix $\mtx{Q}_{\rm orig}$ is called the \textit{update matrix}.
The extended system (\ref{eq:orig_extended_system}) is obtained by subtracting the contributions
from $\Gamma_c$ in $\mtx{A}_{oo}$ and adding the contributions from $\Gamma_p$.  Upon
solving (\ref{eq:orig_extended_system}), only $\vct{\sigma}_k$ and $\vct{\sigma}_p$ are 
used to evaluate the solution inside of $\Gamma_n$.  Effectively $\vct{\sigma}_c$ is a dummy
vector. Details of the derivation of (\ref{eq:orig_extended_system}) are provided in \cite{2009_martinsson_ACTA,ZHANG2018_jcp}.

\subsection{New extended linear system}
\label{sec:new}

\textcolor{dbrown}{The new extended linear system exploits the fact that the contribution from $\Gamma_c$ 
is not used to find the solution inside of $\Gamma_n$.  }
Specifically, we introduce the vector $\vct{\sigma}_c^{\rm dum}$ 
fully knowing a priori that it will contain useless information. Then 
  solving (\ref{equ:perturbedlLinearSystem}) is equivalent to solving the following 
\begin{equation}\label{eq:extendedlLinearSystem}
\begin{bmatrix}
\mtx{A}_{kk} & \mtx{0} & \mtx{A}_{kp}\\
\mtx{A}_{ck} & \mtx{A}_{cc}  & \mtx{0}\\
\mtx{A}_{pk} & \mtx{0} & \mtx{A}_{pp}\\
\end{bmatrix}
\begin{pmatrix*}[l]
\boldsymbol{\sigma}_k\\
\textcolor{black}{\boldsymbol{\sigma}^{\rm dum}_c}\\
\boldsymbol{\sigma}_p
\end{pmatrix*}
=
\begin{pmatrix}
\vct{g}_k\\
\vct{0}\\
\vct{g}_p
\end{pmatrix}\ycmt{.}
\end{equation}

The expanded form of (\ref{eq:extendedlLinearSystem}) is 
\begin{equation}\label{eq:extendedlLinearSystemUpdate}
\left(
\underbrace{
\begin{bmatrix}
\mtx{A}_{oo} & \mtx{0} \\
 \mtx{0} & \mtx{A}_{pp}\\
\end{bmatrix}}_{\tilde{\mtx{A}}} +
\underbrace{
\begin{bmatrix}
 \mtx{0}  & -\mtx{A}_{kc} & \mtx{A}_{kp}\\
 \mtx{0} &  \mtx{0} & \mtx{0}\\
\mtx{A}_{pk} & \mtx{0} & \mtx{0}\\
\end{bmatrix}}_{\mtx{Q}_{\rm new}}
\right)
\underbrace{
\begin{pmatrix*}[l]
\boldsymbol{\sigma}_k\\
\textcolor{black}{\boldsymbol{\sigma}^{\rm dum}_c}\\
\boldsymbol{\sigma}_p
\end{pmatrix*}}_{\textcolor{black}{\boldsymbol{\sigma}_{\rm ext}}}
=
\underbrace{
\begin{pmatrix}
\vct{g}_k\\
\vct{0}\\
\vct{g}_p
\end{pmatrix}}_{\textcolor{black}{\vct{g}_{\rm ext}}}.
\end{equation}
Here $\mtx{Q}_{\rm new}$ is the new update matrix.  Notice that $\mtx{Q}_{\rm new}$ has a zero row.
\textcolor{dblue}{As compared to the original formulation (\ref{eq:orig_extended_system}), the new formulation has two advantages: first, the update matrix no longer contains the full-rank block $\mtx{B}_{cc}$; second the new formulation does not require evaluating $\mtx{A}_{cp}$.
We postpone the discussion for the first advantage to Section \ref{sec:fds}.
For the problems considered in this paper, $N_c$ and $N_p$ are relatively small compared to $N_k$, thus evaluating and compressing $\mtx{A}_{cp}$ is not expensive if the kernel is analytic.}
\textcolor{black}{
However, for problems where the perturbation corresponds to local refinement of the same geometry, the collection of discretization points removed $I_c$ and those added $I_p$ correspond to different discretizations of the same boundary segment $\Gamma_c=\Gamma_p$.  
This means that the evaluation of $\mtx{A}_{cp}$ requires additional care when the kernel is weakly singular.  By not including
this matrix in the new formulation, the problem of handling weakly singular kernels (such as when using a combined field 
integral representation for Helmholtz problems) is removed.
}

\subsection{A fast direct solver}
\label{sec:fds}
When constructing the fast direct solver for the locally perturbed boundary value problem, there are 
  advantages to writing the system in the form of (\ref{eq:orig_extended_system}) and (\ref{eq:extendedlLinearSystemUpdate}).
Since the matrix $\mtx{A}_{oo}$ is the system resulting from the discretization of the integral equation on the 
original geometry, 
\textcolor{dblue}{we assume that a fast direct solver has 
already been computed for $\mtx{A}_{oo}$.}   Any fast direct solver such as 
\textit{Hierarchically Block Separable (HBS)} \cite{2012_martinsson_FDS_survey},
  \textit{Hierarchically Semi-Separable (HSS)}  \cite{2007_shiv_sheng,2004_gu_divide},
   \textit{Hierarchical Interpolative Factorization (HIF)} \cite{HoYing2013} and 
   $\mathcal{H}$ and $\mathcal{H}^2$- matrix methods \cite{hackbusch} can be used.
Additionally, the update matrices $\mtx{Q}_{\rm orig}$ and
$\mtx{Q}_{\rm new}$ are low rank.  This allows for the inverse of the extended systems 
to be applied rapidly via a Sherman-Morrison-Woodbury formula
\textcolor{black}{
\begin{equation}\label{eq:wood}
\boldsymbol{\sigma}_{\rm ext}=
\left( \tilde{\mtx{A}} +\mtx{Q} \right)^{-1}\vct{g}_{\rm ext}\approx
\left( \tilde{\mtx{A}} +\mtx{LR} \right)^{-1}\vct{g}_{\rm ext}
\approx
\tilde{\mtx{A}}^{-1}\vct{g}_{\rm ext}
-\tilde{\mtx{A}}^{-1}\mtx{L}
  \left( \mtx{I} + \mtx{R}\tilde{\mtx{A}}^{-1}\mtx{L}\right )^{-1}
  \mtx{R}\tilde{\mtx{A}}^{-1}\vct{g}_{\rm ext},
\end{equation}}
where $\mtx{I}$ is an identity matrix, and $\mtx{LR}$ denotes the low rank factorization 
of the update matrix $\mtx{Q}$.  

The low rank property of the update matrices $\mtx{Q}_{\rm orig}$ and $\mtx{Q}_{\rm new}$ 
can be observed by noting that the matrices 
$\mtx{A}_{kc}$, $\mtx{A}_{kp}$, $\mtx{A}_{pk}$ and $\mtx{A}_{op}$ are low rank.  The 
only full rank matrix in the update matrices is $\mtx{B}_{cc}$. 
\textcolor{dblue}{
Let $k_{op}$, $k_{kc}$, $k_{pk}$, and $k_{kp}$ be the observed numerical ranks for matrix $\mtx{A}_{op}$, $\mtx{A}_{kc}$, $\mtx{A}_{kp}$, $\mtx{A}_{kp}$ respectively.
Then the low-rank approximation for $\mtx{Q}_{\rm orig}$ has total rank $k_{\rm orig}= k_{op}+k_{kc}+k_{pk}+N_c$.
 Since $\mtx{Q}_{\rm new}$ 
does not contain the $\mtx{B}_{cc}$ matrix, its rank is $k_{\rm new} = k_{kc}+k_{pk}+k_{kp}$, which is smaller than $\mtx{Q}_{\rm orig}$.
Thus the Sherman-Morrison-Woodbury formula can be applied more rapidly with the new formulation. }
The details for 
efficiently creating the low rank factorizations can be found in \cite{ZHANG2018_jcp}.

\section{Numerical experiments}
\label{sec:numerics}

This section illustrates the performance of the fast direct solver for the proposed 
extended linear system for a collection of problems.  
The integral equations are discretized via the Nystr\"om method with a 
\textcolor{black}{$16$-point composite Gaussian 
quadrature}.  For all problems, the original
geometry is discretized with enough points in order for the boundary value problem to 
be solved to 10 digits of accuracy. 
The HBS direct solver \cite{2012_martinsson_FDS_survey} was used in the examples in this 
section.  
For all tests, the tolerance for HBS compression and low-rank approximation is set to $\epsilon=10^{-10}$. 

Roughly speaking the cost of building fast direct solvers is split into two parts: precomputation
and solve.  The time for precomputation is the time for constructing all the parts of the
fast direct solver.  For the extended systems, this includes constructing the low rank 
factorizations of the update matrices $\mtx{Q}_{\rm orig}$ or $\mtx{Q}_{\rm new}$ and 
inverting the small matrix in the Sherman-Morrison-Woodbury formula (\ref{eq:wood}).  For the HBS solver,
the precomputation includes creating a compressed approximation of the discretized 
system on the new geometry and inverting that approximation.  The \textcolor{dbrown}{solve time}
is the time for applying the resulting solver to one vector (or right-hand-side).  For the extended systems, this
is the time for applying the Sherman-Morrison-Woodbury formula (\ref{eq:wood}).  For the HBS 
solver, it is the time for applying the approximate inverse.

To illustrate the efficiency of the proposed technique, we compare the performance of the new solution technique
with the fast solver developed for the original extended system and building an HBS solver
from scratch for the new geometry.  We report the following:
\begin{itemize}
\item $N_o$: the number of discretization points on the original geometry;
\item $N_c$: the number of discretization points cut from the original geometry;
\item \textcolor{black}{$N_k$: the number of discretization points remained the same on the original  and new geometry, $N_k=N_o-N_c$;}
\item $N_p$: the number of discretization points added;
\item \textcolor{dblue}{$k_{op}$, $k_{kc}$, $k_{pk}$, and $k_{kp}$: the observed numerical ranks for matrix $\mtx{A}_{op}$, $\mtx{A}_{kc}$, $\mtx{A}_{kp}$, and $\mtx{A}_{kp}$ respectively after compression;}
\item \textcolor{dblue}{$k_{\rm new}$ and $k_{\rm orig}$: the observed numerical ranks for $\mtx{Q}_{\rm new}$ and $\mtx{Q}_{\rm orig}$ respectively, $\mtx{Q}_{\rm new}= k_{kc} + k_{kp} + k_{pk}$ and $\mtx{Q}_{\rm orig}= k_{op} + k_{kc} + k_{pk} + N_c$ ;}
\item $T_{{\rm new},p}$: the time in seconds for the precomputation of the proposed solver;
\item $T_{{\rm orig},p}$: the time in seconds for the precomputation of the \textcolor{black}{fast solver based on the original extended system formulation in \cite{ZHANG2018_jcp}};
\item $T_{{\rm hbs},p}$: the time in seconds for the precomputation of HBS from scratch for the new geometry;
\item $r_{p}=\frac{T_{{\rm hbs},p}}{T_{{\rm new},p}}$;
\item $T_{{\rm new},s}$: the time in seconds for applying the proposed solver to one right-hand-side;
\item $T_{{\rm orig},s}$: the time in seconds for applying the original solver \textcolor{black}{in \cite{ZHANG2018_jcp}} to one right-hand-side;
\item $T_{{\rm hbs},s}$: the time in seconds for applying the HBS inverse to one right-hand-side;
\item $r_{s}=\frac{T_{{\rm hbs},s}}{T_{{\rm new},s}}$.
\end{itemize}
The ratios $r_p$ and $r_s$ are measures for the speed-up (or slow-down) by using the proposed solver 
versus building a new fast direct solver from scratch for the new geometry.  If $r_p$ is 
greater than 1, the precomputation of the proposed solver is faster than building a fast direct solver
from scratch.  If $r_p$ is less than 1, the precomputation of the proposed solver is 
slower than building a fast direct solver from scratch, etc.

All experiments were run on a dual 2.3 GHz Intel Xeon Processor E5-2695 v3 desktop workstation 
with 256 GB of RAM.  The code 
is implemented in MATLAB, apart from 
the randomized linear algebra utilized in creating low rank factorization rapidly which
is implemented in Fortran.

\subsection{A local change in the geometry}
\label{sec:pino_nose}
\begin{figure}[h]
 \begin{center}
 \begin{picture}(150,120)(100,50)
 \put(20,20){\includegraphics[scale=0.35]{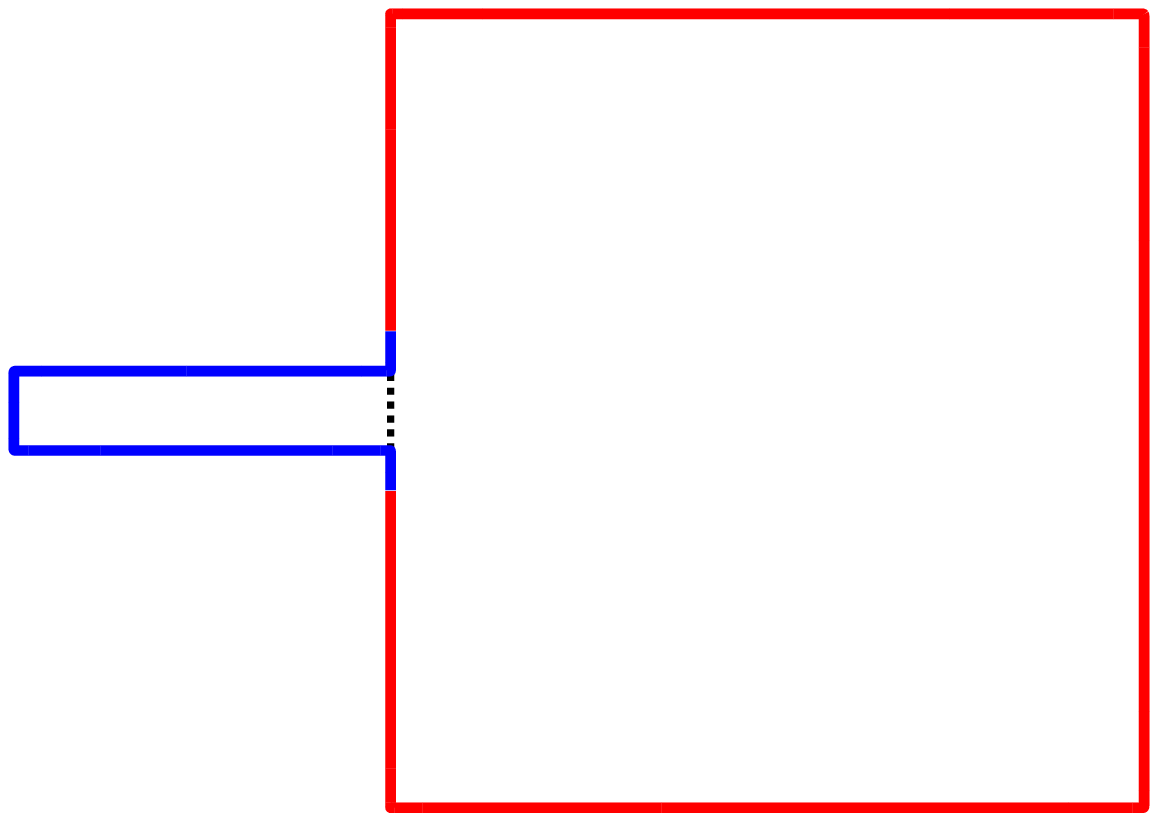}}
 \put(219,70){\includegraphics[scale=0.11]{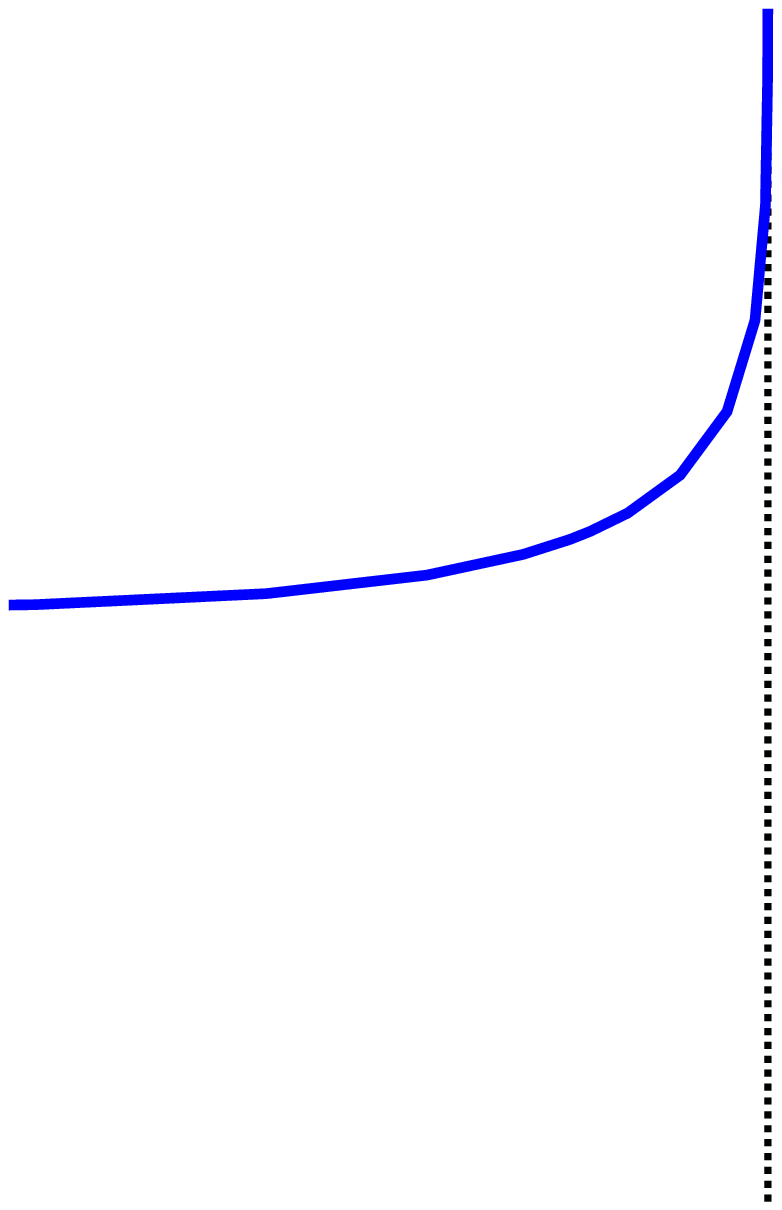}}
 \put(140,72){\footnotesize${\Omega}$}
 \put(116,122){\footnotesize\color{red}${\Gamma}_k$}
 \put(107,90){\footnotesize${\Gamma}_c$}
 \put(240,80){\tiny${\Gamma}_c$}
 \put(235,99){\tiny\color{blue}${\Gamma}_p$}
 \put(68,78){\footnotesize\color{blue}${\Gamma}_p$}
 \put(54,94){\footnotesize$d$}%
 \put(102,100){\circle{6}}
 \qbezier(102,100)(155,120)(220,100)
 \put(220,100){\vector(1,-1){2}}
 \put(245,94){\circle{62}}
 \end{picture}
 \end{center}
 \caption{The square with nose geometry.
 A nose of height $d$ is smoothly attached to the a square.
 }
 \label{fig:squreNoseGeometry}
 \end{figure}

Consider the Laplace boundary value problem (\ref{eq:laplaceBvp}) on the geometry
illustrated in Figure \ref{fig:squreNoseGeometry}. 
The corners are smoothed via the scheme in \cite{EPSTEIN2015_sisc}. A detailed description of this
geometry is given in \cite{ZHANG2018_jcp}. 
\textcolor{black}{The Dirichlet data on the boundary equals to the potential due to
 a collection of $10$ 
charges with location and charge value $\left\{ (\vct{s}_j, \, q_j ) \right\}_{i=1}^{10}$ placed on the exterior of domain $\Omega$, 
$$g(\vct{x}) =\sum_{j=1}^{10} q_j G(\vct{x}, \vct{s}_j),$$
where $G(\vct{x},\vct{y})$ denotes the Green's function for Laplace equation.}

In the first experiment, the number of points cut remains fixed, $N_c = 16$,
while the number of discretization points on $\Gamma_k$ grows.  
In Figure \ref{fig:squreNoseGeometry}, this 
corresponds to the nose height $d$ decreasing as $N_k$ grows.  
\textcolor{dbrown}{The attached nose $\Gamma_p$ is discretized with $N_p\in[832,896]$ number of quadrature points.}
The timing results are reported in Table \ref{tab:thinning_nose_laplace_timing}.  All three solution techniques are linear
with respect to $N_o$ and the precomputation time for the new solution technique is about the same as the original extended system solver.
It is roughly 3.5 times faster than building a new direct solver
from scratch for the new geometry.  
\textcolor{dbrown}{The cost of applying the proposed solver is 
almost as fast as applying the HBS approximate inverse.}
\textcolor{dblue}{To better understand the new extended system formulation's performance as compared to the original formulation  in \cite{ZHANG2018_jcp}, Table \ref{tab:thinning_nose_laplace_rank} reports the numerical ranks for the subblocks  in $\mtx{Q}_{\rm orig}$ and $\mtx{Q}_{\rm new}$. The total rank of $\mtx{Q}_{\rm new}$ is observed to be smaller than that for $\mtx{Q}_{\rm orig}$ by about 20, which is  not a very significant reduction. This explains the timing results in Table \ref{tab:thinning_nose_laplace_timing}, where the timings in column $T_{{\rm orig},p}$ is slightly larger than those in column $T_{{\rm new},p}$ for large tests. }

In the next example,
$N_c$ grows by the same factor as $N_k$. 
\textcolor{dbrown}{The nose height $d$ in Figure \ref{fig:squreNoseGeometry}
remains fixed, 
and the fixed nose is discretized with $N_p=896$ points.}  
Table \ref{tab:fixed_nose_laplace_timing} reports on the 
performance of all three solvers for this geometry.  
The proposed solution
technique is the fastest for the precomputation step.
It is much faster than the solver based on the original extended system formulation, especially for the case where $N_c$ is large.
\textcolor{dblue}{The ranks of the compressed matrices are reported in Table \ref{tab:fixed_nose_laplace_rank}. 
Since $N_c$ grows at the same rate as $N_k$ while the problem size increases, $k_{\rm orig}$ becomes very large 
and the dense inversion of $\left( \mtx{I}+\mtx{R}\tilde{\mtx{A}}^{-1}\mtx{L}\right)$ dominates the cost of evaluating the solution of the original extended system.
This is in constrast to the new formulation where $k_{\rm new}$ does not depend on $N_c$ and thus the new solver is much faster than both the original solver in \cite{ZHANG2018_jcp} and HBS built from scratch for precomputation.}
 \textcolor{dblue}{A factor of roughly 2.9 speed up in the precomputation is observed as compared to an HBS solver built from scratch.}  Again applying the proposed
solver is slightly slower than applying the HBS approximate inverse. 

%
%
 \begin{table}[h]
 \begin{center}
 \begin{tabular}{lllllllll}
 \hline
 $N_o$ &  $T_{{\rm orig},p}$ &  $T_{{\rm new},p}$ & $T_{{\rm hbs},p}$ & $r_p$
    & $T_{{\rm orig},s}$& $T_{{\rm new},s}$ & $T_{{\rm hbs},s}$ & $r_s$   \\ \hline
 9232            
 & 3.69e-01		& 4.83e-01 & 1.57e+00  & 3.25 
 & 1.99e-02 	& 1.12e-02 & 1.32e-02  & 1.18 \\ 
 18448          
 & 5.60e-01 		& 6.50e-01 & 2.38e+00  & 3.66 
 & 2.76e-02  		& 1.74e-02 & 1.46e-02 & 0.84 \\ 
 36880         
 &  1.11e+00      	& 1.11e+00 & 3.79e+00 & 3.42 
 & 5.49e-02			& 4.00e-02 & 3.33e-02 & 0.83 \\
 73744          
 & 2.25e+00 		& 1.84e+00 & 6.38e+00  & 3.47 
 & 9.79e-02 		& 8.06e-02 & 7.04e-02 & 0.87 \\ 
 147472        
 & 3.87e+00    		& 3.56e+00 & 1.18e+01 & 3.33 
 & 1.95e-01  		& 1.71e-01 & 1.52e-01 & 0.89 \\ 
 
 \hline
 \end{tabular}
 \caption{Times for applying the solution technique to (\ref{eq:laplaceBvp}) on the square with thinning nose geometry.
  }
 \label{tab:thinning_nose_laplace_timing}
 \end{center}
 
 \end{table}

 \begin{table}[h]
 \begin{center}
 \begin{tabular}{llllll|l|l}
 \hline
 $N_o$ & $N_c$ & $k_{op}$ &  $k_{kc}$ & $k_{pk}$ & $k_{kp}$
    & $k_{\rm orig}$& $k_{\rm new}$  \\ \hline
 9232 &  16          
 & 63 & 9 & 54 & 60 & 142 & 123\\ 
 18448 &  16         
 & 64 & 10 & 58 & 58 & 148 & 126\\ 
 36880  &  16       
 &  66 & 10 & 57 & 61 & 149 & 128\\
 73744  &  16        
 &  68 & 10 & 59 & 60 & 153 & 129\\ 
 147472  &  16      
 &  69 & 10 & 61 & 61 & 156 & 132\\ 
 
 \hline
 \end{tabular}
 \caption{ \textcolor{dblue}{
Observed numerical ranks for compressed matrices in the original and new extended system formulations for the square with thinning nose geometry.
  }}
 \label{tab:thinning_nose_laplace_rank}
 \end{center}
 \end{table}

 \begin{table}[h]
 \begin{center}
 \begin{tabular}{llllllllll}
 \hline
 $N_o$ & $N_c$  &  $T_{{\rm orig},p}$&$T_{{\rm new},p}$ & $T_{{\rm hbs},p}$ & $r_p$
      & $T_{{\rm orig},s}$& $T_{{\rm new},s}$ & $T_{{\rm hbs},s}$ & $r_s$    \\ \hline
 9344     &   128 
 &	5.01e-01    & 5.10e-01 & 1.28e+00 & 2.50 
 &	2.08e-02	& 1.02e-02 & 7.92e-03&  0.77 \\ 
 18688    &   256
 &	1.08e+00        & 9.25e-01 & 2.18e+00 & 2.36
 &	3.44e-02   		& 2.15e-02 & 1.59e-02 & 0.74 \\ 
 37376    &   512 
 &2.67e+00      	& 1.30e+00 & 3.49e+00 & 2.69
 &5.64e-02   		&  3.97e-02 & 3.00e-02 & 0.76 \\ 
 74752    &   1024  
 &7.76e+00   		& 2.31e+00 & 6.63e+00 & 2.87
 &1.16e-01 		& 8.67e-02 & 6.40e-02 & 0.74 \\
 149504  &   2048  
 &2.48e+01      	&  4.06e+00 & 1.19e+01  & 2.92 
 &2.34e-01	 		& 1.71e-01 & 1.61e-01& 0.94\\ 
  \hline
 \end{tabular}
 \caption{Times for applying the solution techniques to (\ref{eq:laplaceBvp}) on the square with fixed nose geometry.
}
 \label{tab:fixed_nose_laplace_timing}
 \end{center}
 
 \end{table}
 
\begin{table}[h]
 \begin{center}
 \begin{tabular}{llllll|l|l}
 \hline
 $N_o$ & $N_c$ & $k_{op}$ &  $k_{kc}$ & $k_{pk}$ & $k_{kp}$
    & $k_{\rm orig}$& $k_{\rm new}$  \\ \hline
 9344     &   128 
 &	87 & 12 & 45 & 53 & 272 & 99\\ 
 18688    &   256
 &	80 & 14 & 47 & 55 & 397 & 116\\ 
 37376    &   512 
 & 96 & 14 & 45 & 55 & 667 & 114\\ 
 74752    &   1024  
 & 117 & 16 & 46 & 59 & 1203 & 121\\
 149504  &   2048  
 & 125 & 18 & 46 & 56 & 2237 &120\\ 
 \hline
 \end{tabular}
 \caption{\textcolor{dblue}{Observed numerical ranks for compressed matrices in the original and new extended system formulations for the square with fixed nose geometry.  }}
 \label{tab:fixed_nose_laplace_rank}
 \end{center}
 
 \end{table}

\subsection{A Laplace problem with a locally refined discretization}
\begin{figure}[h]
\begin{center}
\begin{picture}(300,170)(01,01)
\put(-100,10){\includegraphics[scale=0.45]{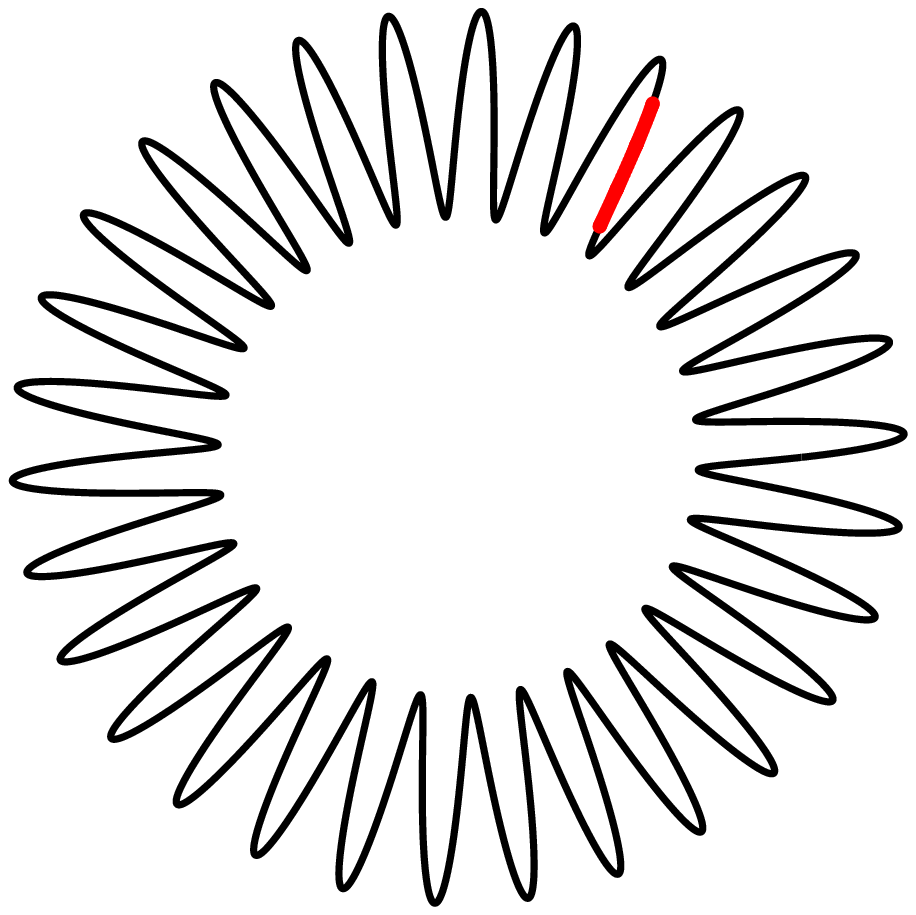}}
\put(45,40){$\Gamma_k$}
\put(-10,10){(a)}
\put(80,21){\includegraphics[scale=0.35]{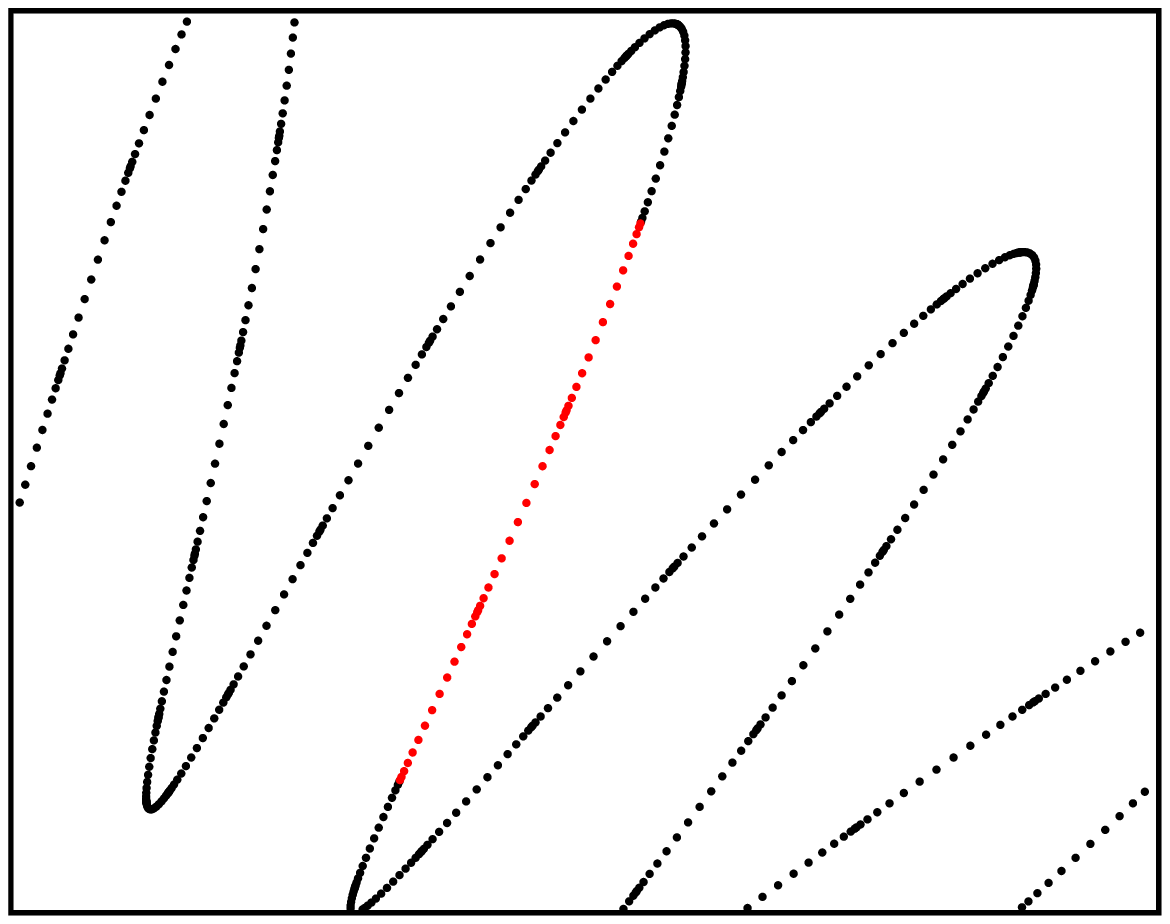}}
\put(140,10){\footnotesize (b)}
\put(220,21){\includegraphics[scale=0.35]{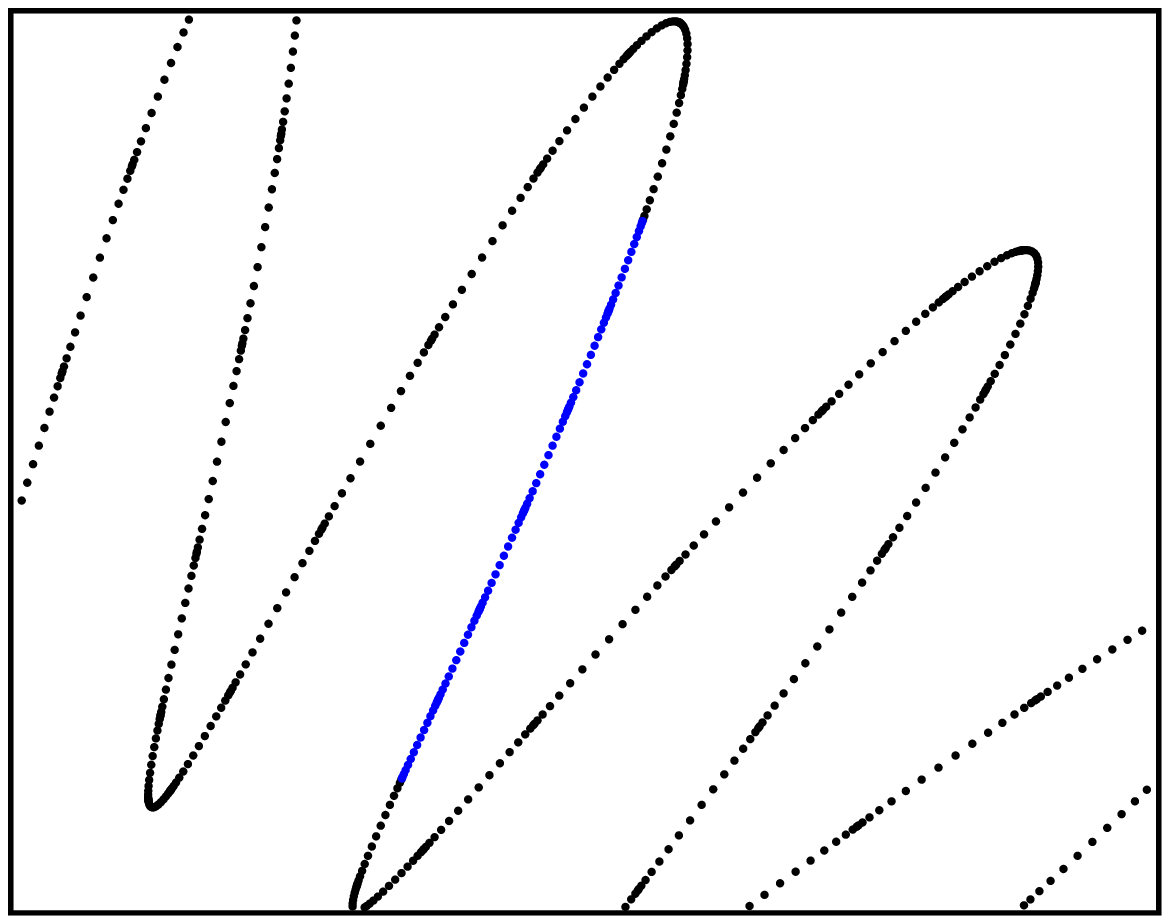}}
\put(280,10){\footnotesize (c)}
\put(140,80){\color{red}$\Gamma_c$}
\put(280,80){\color{blue}$\Gamma_p$}

\end{picture}
\end{center}
\caption{(a) The sunflower geometry with the portion of the boundary to 
be refined in red.  (b) The three Gaussian panels in the boxed region
from the original discretization. (c)  Six Gaussian panels replacing
the original three panels.
}
\label{fig:sunflower_geom}
\end{figure}

Next we consider applying the proposed solution technique to the Laplace boundary 
value problem (\ref{eq:laplaceBvp}) where 
the local perturbation is a refinement in a portion of the geometry.  Figure
\ref{fig:sunflower_geom}(a) illustrates the geometry under consideration.  It is
given by the following parameterization:
$$\mtx{x}(t)=
\begin{pmatrix}
r(t) \cos(t)\\
r(t) \sin(t)
\end{pmatrix}
,\mbox{ with }
r(t)=1+0.3\sin(30t)
\mbox{ for }t\in[0,2\pi].
$$
The portion of the boundary being refined is highlighted in red. Figure \ref{fig:sunflower_geom}(b)
is a zoomed in illustration of that region. Figure \ref{fig:sunflower_geom}(c) illustrates
the local refinement.  Three Gaussian panels ($N_c=48$) are replaced with $N_p$ discretization points
($N_p/16$ Gaussian panels). The number of discretization points on $\Gamma_k$ remains fixed; $N_k = 6352$.
\textcolor{black}{The Dirichlet data is generated similarly as in Section \ref{sec:pino_nose}.}

Table \ref{tab:sunflower_laplace_timing} reports on the performance of all three solution 
techniques for this problem.  The proposed solution technique is $13$ to $21$ times faster
than building a new solver from scratch while applying the solver is less than a factor two
slower than applying the HBS approximate inverse.

\begin{table}[h]
\begin{center}
\begin{tabular}{llllllllll}
\hline
$N_p$ & $\frac{N_p}{N_o}$ & $T_{{\rm orig},p}$ & $T_{{\rm new},p}$ & $T_{{\rm hbs},p}$ & $r_p$   & $T_{{\rm orig},s}$  & $T_{{\rm new},s}$ & $T_{{\rm hbs},s}$ & $r_s$    \\ \hline
96    & 0.015         & 6.06e-01  	    & 5.03e-01 & 7.52e+00 & 14.9  &1.10e-02   & 1.30e-02 & 1.32e-02 & 1.02 \\ 
192   & 0.03           &6.16e-01	       & 3.62e-01 & 7.77e+00 & 21.4  &1.17e-02	   & 1.25e-02 & 9.30e-03 & 0.74  \\ 
384   & 0.06            &6.83e-01	      & 3.90e-01 & 7.72e+00 & 19.8   &1.36e-02	  & 1.42e-02 & 9.13e-03 & 0.64 \\ 
768   & 0.12            & 7.60e-01	     & 4.11e-01 & 7.78e+00 & 18.9   &2.01e-02	  & 1.20e-02 & 9.06e-03 & 0.76 \\ 
1536  & 0.24            &	1.01e+00      & 6.09e-01 & 8.03e+00 & 13.2   &4.72e-02	  & 1.66e-02 & 1.00e-02 & 0.60 \\ 
\hline
\end{tabular}
\caption{Times for applying the solution techniques to (\ref{eq:laplaceBvp}) on the geometry in 
Figure \ref{fig:sunflower_geom} with local refinement.}
\label{tab:sunflower_laplace_timing}
\end{center}

\end{table}

\subsection{A Helmholtz problem with a locally refined discretization}
\label{sec:helmholtz}
Besides being faster than the solver for the original extended system, the proposed
solver has the advantage that it can easily handle problems that are using 
specialized quadrature for weakly singular kernels.  The issue that arises 
for the original extended system is that it would be cumbersome to evaluate 
the entries of the matrix $\mtx{A}_{op}$ corresponding to the interaction
of $\Gamma_c$ with $\Gamma_p$.  This matrix does not arise in the 
new extended system. 

To illustrate the efficiency of the solver for systems that involve specialized 
quadrature, we consider the following exterior Dirichlet Helmholtz boundary value 
problem
\begin{equation}\label{eq:helmholtz}
\begin{array}{rlll}
 -\Delta u(\vct{x}) +\omega^2u &=&0 &\mbox{ for }\vct{x} \in \Omega^c,\\
 u(\vct{x})&=&g(\vct{x})&\mbox{ for }\vct{x} \in \Gamma\\
\end{array}
\end{equation}
 with Sommerfeld radiation condition on the sunflower geometry illustrated
in Figure \ref{fig:sunflower_geom} where $\omega$ denotes the wave number.
\textcolor{black}{
The Dirichlet data $g$ is set to be the negative of a plane wave with incident angle $\theta=-\frac{\pi}{5}$
$$g(\vct{x})=-e^{i\vct{k}\cdot \vct{x}},\mbox{ with }\vct{k}=(\omega \cos\theta, \omega \sin\theta).$$}

We chose to represent the solution with the following combined field 
\begin{equation}
 u(\vct{x})=\int_\Gamma D_\omega(\vct{x},\vct{y}) \sigma(\vct{y})\, ds(\vct{y}) - i \omega \int_\Gamma S_\omega(\vct{x},\vct{y}) \sigma(\vct{y})\, ds(\vct{y}),
\end{equation}
\textcolor{dbrown}{
 where $\sigma(\vct{x})$ is the unknown boundary charge distribution, $S_\omega=G_{\omega}(\vct{x},\vct{y})$ and 
  $ D_\omega=\partial_{\vct{n}({\vct{y}})} G_\omega(\vct{x},\vct{y})$ denote the single and double layer Helmholtz kernel, 
  $\vct{n}({\vct{y}})$ is the outward facing normal vector,
and  $G_{\omega}(\vct{x},\vct{y}) = \frac{i}{4}H_0^{(1)}(\omega|\vct{x}-\vct{y}|)$ is the two dimensional free space Green's function for the Helmholtz equation with wave number $\omega$ and $H_0^{(1)}$ is the Hankel function of zeroth order.
}

The integral equation that results from enforcing the Dirichlet boundary condition is 
\begin{equation}
\frac{1}{2} \sigma(\vct{x})+ \int_\Gamma D_\omega(\vct{x},\vct{y}) \sigma(\vct{y})\, ds(\vct{y}) - i \omega \int_\Gamma S_\omega(\vct{x},\vct{y}) \sigma(\vct{y})\, ds(\vct{y})=g(\vct{x}).
\end{equation}

We discretize the operator via Nystr\"om with a composite
\textcolor{black}{$16$-point generalized Gaussian quadrature} 
\cite{2001_rokhlin_kolm}.  The wave number is set to $\omega=20$ which corresponds to the geometry being
approximately 8.3 wavelengths in size. 
 Again, we consider the local refinement problem.
Table \ref{tab:sunflower_helmholtz_timing} reports on the performance of the proposed
solution technique and building a fast direct solver from scratch.  For this problem,
the proposed solver is anywhere from 15 to 35 times faster than building the fast direct
solver from scratch.  This speed up is the result of the increased ranks associated with 
Helmholtz problems.  Applying the proposed solver to a right-hand-side is roughly
1.5 times slower than applying the HBS solver.

\begin{table}[h]
\begin{center}
\begin{tabular}{llllllll}
\hline
$N_p$ & $\frac{N_p}{N_o}$ & $T_{{\rm new},p}$ & $T_{{\rm hbs},p}$ & $r_p$ & $T_{{\rm new},s}$ & $T_{{\rm hbs},s}$ &$r_s$    \\ \hline
96    & 0.015             & 1.13e+00 & 3.97e+01 & 35.2 & 4.06e-02 & 2.86e-02 & 0.71 \\ 
192   & 0.03              & 1.36e+00 & 4.08e+01 & 29.9 & 4.64e-02 & 2.93e-02 & 0.63 \\ 
384   & 0.06              & 1.44e+00 & 4.08e+01 & 28.4 & 3.91e-02 & 2.54e-02 & 0.65\\ 
768   & 0.12               & 1.64e+00 & 4.17e+01 & 25.4 & 3.69e-02 & 2.70e-02 & 0.73 \\ 
1536  & 0.24            & 2.64e+00 & 4.08e+01 & 15.4 & 4.39e-02 & 3.33e-02 & 0.76 \\
 \hline
\end{tabular}
\caption{Times for applying the solution techniques to (\ref{eq:helmholtz}) on the geometry in 
Figure \ref{fig:sunflower_geom} with local refinement.}
\label{tab:sunflower_helmholtz_timing}
\end{center}
\end{table}

\section{Concluding remarks}  This manuscript presented a new extended linear system for integral equation based 
solution techniques for boundary value 
problems on locally perturbed geometries.  
\textcolor{black}{A fast direct solver based on the new extended system formulation is 
significantly faster than building a new direct solver from scratch for the perturbed problem.
For some examples, the precomputation is between 10x to 30x faster than building a direct 
solver from scratch. Additionally, the new solver shows consistent speed-ups for problems with a large number of points removed and can be easily
applied to problems requiring the discretization of weakly singular kernels.  Neither of these was true 
for the fast direct solver based on the original formulation given in \cite{ZHANG2018_jcp}. }

\textcolor{dblue}{The idea of handling local changes in geometry or discretization via the extended system formulation can be extended to three dimensional problems. However, additional work is required in  processing the geometry and 
creating efficient techniques for building the low rank factors of the update matrix.  This is future work.
   }


\begin{thebibliography}{10}

\bibitem{2004_gu_divide}
S.~Chandrasekaran and M. Gu
\newblock A divide-and-conquer algorithm for the eigendecomposition of symmetric block-diagonal plus semiseparable matrices.
\newblock {\em Numerische Mathematik}, 96(4):723-731, 2004

\bibitem{EPSTEIN2015_sisc}
C. ~Epstein and M. ~O'Neil.
\newblock  Smoothed Corners and Scattered Waves.
\newblock {\em SIAM Journal on Scientific Computing}, 38, 2015


\bibitem{2012_martinsson_FDS_survey}
A.~Gillman, P.~Young and P.G.~Martinsson
\newblock A direct solver $O(N)$ complexity for integral equations on one-dimensional domains.
\newblock {\em Frontiers of Mathematics in China}, 7:217-247, 2012

\bibitem{2009_martinsson_ACTA}
L.~Greengard, D.~Gueyffier, P.G.~Martinsson and V. Rokhlin
\newblock Fast direct solvers for integral equations in complex three-dimensional domains.
\newblock {\em Acta Numerica}, 18:243-275, 2009


\bibitem{hackbusch}
W.~Hackbusch
\newblock A Sparse Matrix Arithmetic Based on {H}-Matrices; {P}art {I}: {I}ntroduction to {H}-Matrices.
\newblock {\em Computing}, 62:89-108, 1999

\bibitem{HoYing2013}
K.~Ho and L.~Ying
\newblock Hierarchical Interpolative Factorization for Elliptic Operators: Integral Equations.
\newblock {\em Communications on Pure and Applied Mathematics}, 69, 2013





\bibitem{2001_rokhlin_kolm}
P.~Kolm and V. Rokhlin
\newblock Numerical quadratures for singular and hypersingular integrals.
\newblock {\em Computers and Mathematics with Applications}, 41:327-352, 2001



\bibitem{2007_shiv_sheng}
Z. ~Sheng, P. ~Dewilde and S. ~Chandrasekaran
\newblock Algorithms to solve hierarchically semi-separable systems.
\newblock {\em Operator Theory: Advances and Applications }, 176: 255-294, 2007

\bibitem{ZHANG2018_jcp}
Y.~Zhang and A.~Gillman.
\newblock A fast direct solver for boundary value problems on locally perturbed geometries.
\newblock {\em Journal of Computational Physics}, 356: 356 - 371, 2018








\end{thebibliography}
\end{document}